\documentclass[11pt]{amsart}

\usepackage[english]{babel}
\usepackage[latin1]{inputenc}
\usepackage{amsthm,amsmath,amsfonts, amssymb}

\def \ZZ{\mathbb Z}
\def \RR{\mathbb R}
\def \EE{\mathbb E}
\def \PP{\mathbb P}

\def \QQ{\mathbb Q}

\def \C{\mathcal C}
\def \F{\mathcal F}
\def \S{\mathcal S}
\def \N{\mathcal N}
\def \L{\mathcal L}

\def \indic#1{1\!\!1_{\!#1}}
\def \open#1{#1^{o}}
\def \close#1{\overline{#1}}

\theoremstyle{plain}
\newtheorem{theo}{Theorem}[section]
\newtheorem{prop}[theo]{Proposition}
\newtheorem{lemm}[theo]{Lemma}

\theoremstyle{definition}
\newtheorem{rema}[theo]{Remark}
\newtheorem{exem}[theo]{Example}

\begin{document}

\author[R.~Garbit]{Rodolphe Garbit}
\address{Universit\'e d'Angers\\D\'epartement de Math\'ematiques\\ LAREMA\\ UMR CNRS 6093\\ 2 Boulevard Lavoisier\\49045 Angers Cedex 1\\ France}
\email{rodolphe.garbit@univ-angers.fr}

\title[Two-dimensional random walks in a cone]{A central limit theorem for two-dimensional random walks  in a cone}

\thanks{This work was partially supported by Agence Nationale de la Recherche Grant ANR-09-BLAN-0084601}

\date{\today}

\begin{abstract}
We prove that a planar random walk with bounded increments and mean zero which is conditioned to stay in a cone converges weakly to the corresponding Brownian meander if and only if the tail distribution of the exit time from the cone is regularly varying. This condition is satisfied in many natural examples.
\end{abstract}

\maketitle

\section{Introduction}
\subsection{Main result} 
The aim of this paper is to underscore a natural necessary and sufficient condition for the weak convergence of a two-dimensional random walk conditioned to stay in a cone to the corresponding Brownian meander. The condition only involves the asymptotic behavior of the tail distribution of the first exit time from the cone.

Let $(\xi_n)_{n\geq 1}$ be a sequence of independent and identically distributed random vectors of $\RR^d$, $d\geq 1$, defined on a probability space $(\Omega, \F,\PP)$. 
We assume that the distribution of $\xi_1$ satisfies
$\EE(\xi_1)=0$ and $Cov(\xi_1)=\sigma^2I_d$, where $\sigma^2>0$ and $I_d$ is the $d\times d$ identity matrix.

We form the random walk $S=(S_n)_{n \geq 1}$ by setting $S_n=\xi_1+\cdots+\xi_n$, and
for each $n\geq 1$, we define a normed and linearly interpolated version of $S$ by
$$\S_n(t)=\frac{S_{[nt]}}{\sigma\sqrt{n}}+(nt-[nt])\frac{\xi_{[nt]+1}}{\sigma\sqrt{n}},\quad t\geq 0,$$
where $[a]$ denotes the integer part of $a$.

The weak convergence of the process $\S_n = (\S_n(t), t \geq 0)$ as $n\to\infty$ to a standard
Brownian motion is Donsker's theorem (see for example Theorem 10.1
of~\cite{Bil99}).

We consider a linear cone $C\subset\RR^d$  ({\em i.e.} $\lambda C=C$ for every $\lambda>0$) with the following properties:
\begin{enumerate}
\item $C$ is convex,
\item its interior $\open{C}$ is non-empty,
\item $\PP(\xi_1\in C\setminus\{0\})>0$.
\end{enumerate}
Such a cone is said to be {\em adapted} to the random walk. Note that the convexity of $C$ ensures that its boundary $\partial C$ is negligible with respect to Lebesgue measure (see for example~\cite{Lan86}).
The third condition ensures that the first step of the random walk is in $C$ with positive probability. Since a convex cone is a semi-group, the event $\{\xi_1,\ldots,\xi_n\in C\}$ is a subset of $\{S_1,\ldots,S_n\in C\}$, so the latter has also a positive probability. For this purpose, one could simply require that $\PP(\xi_1\in C)>0$, but our third condition also excludes the uninteresting cases where 
$\{S_1,S_2,\ldots, S_n\in C\}=\{S_1=S_2=\cdots=S_n=0\}$ almost surely.

We consider the first exit time of the random walk from the cone defined by
$$T_C=\inf\{n\geq 1:S_n\notin C\},$$
and wish to investigate the asymptotic distribution of $(S_1,\ldots,S_n)$ conditional on $\{T_C>n\}$ as $n\to\infty$. 

We denote by $\C_1$ the space of all continuous functions $w:[0,1]\to \RR^d$, endowed with the topology of the uniform convergence and the corresponding Borel $\sigma$-algebra. Weak convergence of probability measures on $\C_1$ will be denoted by the symbol $\Rightarrow$.

Let $Q_n$ denote the distribution on $\C_1$ of the process $\S_n$ conditional on $\{T_C>n\}$, that is, for any Borel set $B$ of $\C_1$,
$$Q_n(B)=\PP(\S_n\in B\vert T_C>n).$$
Note that, since $C$ is a convex cone, this is equivalent to conditioning $\S_n$ on $\{\tau_C(\S_n)>1\}$, where
$$\tau_C(w)=\inf\{t>0: w(t)\notin C\},\quad w\in \C_1.$$

We are interested in the weak convergence of the sequence of conditional distributions $(Q_n)$.
The one-dimensional case, where $C=[0,\infty)$, has been investigated in the 60's and the 70's by many authors.
It was Spitzer~\cite{Spi60} who first announced a central limit theorem for the random walk conditioned to stay positive:
$$Q_n(w(1)\leq x)\to 1-\exp(-x^2/2),\quad x\geq 0.$$
But, apparently, he never published the proof. Note that the limit is the Rayleigh distribution.
A first proof of the weak convergence of $Q_n$ was given by Iglehart in~\cite{Igl74} under the assumptions $\EE(\vert \xi_i\vert^3)<\infty$ and $\xi_i$ nonlattice or integer valued with span~$1$. The limit is found to be the distribution of Brownian meander. Then Bolthausen proved in~\cite{Bol76} that these extra assumptions were superfluous.
For the reader who may not be familiar with the Brownian meander, we will use a theorem of Durrett, Iglehart and Miller~\cite{DIM77} as a definition.
Let $W^x$ be the distribution of the standard Brownian motion started at $x$. For any $x>0$, we denote by $M^x$ the distribution $W^x$ conditional on $\{\tau_C>1\}$, that is  
$$M^x(B)=W^x(B\vert \tau_C>1)$$
for any Borel set $B$ of $\C_1$.
Here, the definition of conditional probabilities is elementary since $W^x(\tau_C>1)$ is positive for all $x>0$.
The distribution $M$ of Brownian meander is the weak limit of $M^x$ as $x\to 0^+$ (see~\cite{DIM77}, Theorem 2.1). 
Note that the existence of a limit is not straightforward since $W^0(\tau_C>1)=0$.
But in a sense, the Brownian meander is a Brownian motion started at $0$ and conditioned to stay positive for a unit of time. The Brownian meander can alternatively be obtained by some path transformations of Brownian motion. Namely, it is the first positive excursion of Brownian motion with a lifetime greater than $1$; it is also the absolute value of the rescaled section of Brownian motion observed on the interval $[h,1]$, where $h$ is its last zero before $t=1$ (see~\cite{Bol76} and~\cite{DIM77}).

With this in mind, the weak convergence of $Q_n$ to $M$ can be stated in the following imprecise but intuitive way: the random walk conditioned to stay positive converges to a Brownian motion conditioned to stay positive.

We now turn to the two-dimensional case. If $Q_n$ does converge weakly, then its limit should naturally be the distribution of a Brownian motion conditioned to stay in the cone $C$ for a unit of time. Such a process can be defined as the weak limit of conditioned Brownian motion in the same way as Brownian meander. As above, for $x\in\open{C}$, let $M^x$ be the distribution of Brownian motion started at $x$ and conditioned to stay in $C$ for a unit of time. The following theorem is due to Shimura~\cite{Shi85} and has been extended in~\cite{Gar09} to any dimension $d\geq 2$ for smooth cones.

\begin{theo}[\cite{Shi85}, Theorem 2]\label{convmeander} 
As $x\in\open{C}\to 0$, the distribution $M^x$ converges weakly to a limit $M$.
\end{theo}
The limit distribution $M$ in this theorem will be referred to as the distribution of the Brownian meander (of the cone $C$). We will give more details about $M$ in Section~\ref{preparatory}.

We now come to the main result of the present paper. We recall that a sequence $(u_n)$ of positive numbers is {\em regularly varying} if it can be written as $u_n=n^{-\alpha}l_n$, where $\alpha\in\RR$  and $(l_n)$ is  slowly varying, {\em i.e.} $\lim_n l_{[nt]}/l_n= 1$ for all $t>0$ (see for example~\cite{Boj73}). The exponent $\alpha$ is unique and called the index of regular variation.
A non-increasing sequence of positive numbers $(u_n)$ will be called {\em dominatedly varying}\footnote{This is strictly weaker than regular variation. For example, since $\prod_n (1+1/n)$ is divergent, it is possible to construct a sequence of numbers $1\leq c_n\leq 2$ such that :(i) for all $n\geq 1$,
$c_{n+1}\leq (1+1/n)c_n$, and (ii) for all $\epsilon>0$, there exist infinitely many $n$ such that $c_n\geq 2-\epsilon$ and $c_{n+1}=1$. Then, the sequence $u_n=c_{n}/n$ is non-increasing and dominatedly varying, but not regularly varying since $\liminf u_{n+1}/u_n\leq 1/2 $ is not equal to $1$ as it should be (see~\cite{Wei76} for a very nice proof of this).}
if $\limsup_n u_{[nt]}/u_n$ is finite for all $t\in(0,1]$.

\begin{theo}\label{mainthm}
Assume that the two-dimensional random walk has bounded increments. Then, the
sequence of conditional distributions $(Q_n)$ converges weakly to the Brownian meander if and only if $\PP(T_C>n)$ is dominatedly varying.
In that case, $\PP(T_C>n)$ is regularly varying with index $\pi/(2\beta)$, where $\beta$ is the angle of the cone.
\end{theo}

The assumption of bounded increments is only used in the proof of the tightness of $(Q_n)$ which is taken from the paper~\cite{Shi91} of Shimura.
We will discuss some extensions to the case where the increments are not bounded in Section~\ref{approach}.
However, the rest of the proof of Theorem~\ref{mainthm}, which consists in a study of the (eventual) limit points of the sequence $(Q_n)$, is completely independant of the assumption of bounded increments. Thus, we could have stated a more general (but not very useful) theorem by simply assuming that $(Q_n)$ is tight.
In order to avoid any confusion, the reader is advised that in any of the lemmas, propositions or theorems of this paper, the random walk $(S_n)$ is not assumed to have bounded increments unless it is written explicitly. 

Our Theorem~\ref{mainthm} can be regarded as an extension of a previous result due to Shimura (\cite{Shi91}, Theorem~1). Indeed, he proved that $Q_n\Rightarrow M$ if the distribution of the increments satisfies the following condition: there exists an orthogonal basis $\{\vec{u},\vec{v}\}$ of $\RR^2$ with $\vec{v}\in\open{C}$ such that $\EE(V\vert U)=0$, where $(U,V)$ denotes the coordinates of $\xi_1$ in the new basis.
But this condition does not seem to be very natural. For example, consider the simple random walk on $\ZZ^2$. It is well known that the coordinates $(U,V)$ of $\xi_1$ in the basis $\{\vec{u},\vec{v}\}=\{(1,-1),(1,1)\}$ are independent, therefore Shimura's theorem applies to any of the cones
$\{(x,y):0\leq y\leq rx\}$, $r>1$. However, if $C$ is the octant $\{(x,y):0\leq y\leq x\}$, there is no $\vec{v}$ in $\open{C}$ which satisfies the assumption of his theorem.

By comparison, regarding the example of the simple random walk, our Theorem~\ref{mainthm} combined with the precise estimates of $\PP(T_C>n)$ given by Varopoulos in~\cite{Var99} shows that the weak convergence $Q_n\Rightarrow M$ holds for every adapted cone (Theorem~\ref{SR}). Indeed, Varopoulos estimates enable us to state an invariance principle (Theorem~\ref{PI}) that holds for a large class of random walks and adapted cones. 

\subsection{Examples}\label{examples}

Let $\L$ be the set of probability measures $\mu$ on $\RR^2$ with  bounded support, mean zero, covariance matrix $\sigma^2I_2$ with $\sigma^2>0$, and such that either $\mu$ has its support on $\ZZ^2$, or $\mu$ is absolutely continuous with respect to Lebesgue measure. Assume that the distribution $\mu$ of $\xi_1$ belongs to $\L$. Recall $S=(S_n)_{n\geq 1}$ is the associated random walk and let $supp S$ denote its support.

If $supp S\subset \ZZ^2$, define $p(n,x,y)$ by
$$p(n,x,y)=\PP^x(S_n=y;T_C>n),\quad x,y\in\ZZ^2,$$ 
where $\PP^x(S\in *)$ stands for $\PP(x+S\in *)$ as usual.

If $\mu$ is absolutely continuous, define $p(n,x,y)$ by
$$p(n,x,y)dy=\PP^x(S_n\in dy;T_C>n),\quad x,y\in\RR^2,$$
that is $y\mapsto p(n,x,y)$ is the density of the measure $\PP^x(S_n\in *;T_C>n)$.

Following the terminology of~\cite{Var99}, we shall say that an adapted cone $C$ is in {\em general position} with respect to $\mu$ if
\begin{enumerate}
\item[A)] for all $a>0$, there exist $n\geq 1$ and $\epsilon>0$ such that
$$\forall x,y\in C\cap supp S,\quad  \Vert x-y\Vert\leq a\quad\mbox{ implies }\quad p(n,x,y) \geq\epsilon.$$
\item[B)] for all $a>0$, there exist $n\geq 1$ and $\epsilon>0$ such that
$$\forall x\in C\cap supp S,\quad  d(x,\partial C)\leq a\quad\mbox{ implies }\quad \PP^x(T_C\leq n) \geq\epsilon.$$
\end{enumerate}

Assuming that $C$ is in general position with respect to $\mu$, Varopoulos obtained in~\cite{Var99} precise estimates of the tail distribution of $T_C$. In particular, his results show that there exists $\gamma>1$ such that
$$\gamma^{-1}n^{-\pi/2\beta}\leq\PP(T_C>n)\leq\gamma n^{-\pi/2\beta},$$
where $\beta$ is the angle of the cone. This, combined with Theorem~\ref{mainthm}, proves the following invariance principle.

\begin{theo}\label{PI}
Assume that the distribution of $\xi_1$ belongs to $\L$, and that $C$ is in general position. Then $(Q_n)$ converges weakly to $M$.
\end{theo}

The reader may have noticed that the simple random walk on $\ZZ^2$ does never fit the condition A) above because of a parity problem: if $y$ and $z$ are two neighbours, then $p(n,x,y)$ and $p(n,x,z)$ can not be positive simultaneously. However, it is possible to get around this and find
a good estimate of $\PP(T_C>n)$ in this context. Let us explain how.

Let $(S_n)$ be the simple random walk on $\ZZ^2$ and let $C$ be an adapted cone.
For any $x\in\ZZ^2\cap C$, one can always find a path (on $\ZZ^2$) that joins $0$ to $x$ without leaving $C$, hence there exists $k\geq 0$ such that $\PP(S_k=x,T_C>k)>0$.
Fix $x\in\ZZ^2$ deep inside the cone so that $d(x+C,\partial C)>1$. Since the steps of $(S_n)$ are bounded by $1$, we have
$$\{S_{2i}\in C\mbox{ for all }i=1\ldots m\}\subset\{S_j\in C-x\mbox{ for all }j=1\ldots 2m\},$$
and consequently
$$\PP(S_{2i}\in C\mbox{ for all }i=1\ldots m)\leq\PP^{x}(S_j\in C\mbox{ for all }j=1\ldots 2m).$$
Now, choose $k$ such that $\alpha:=\PP(S_k=x,T_C>k)>0$ and set $m=[\frac{n-k}{2}]+1$.
Using the Markov property of the random walk and the preceding inequality, we obtain, for $n>k$,
\begin{align}\label{up}
\PP(T_C>n)&\geq \PP(S_k=x,T_C>n)\nonumber\\
&\geq \PP(S_k=x,T_C>k)\PP^{x}(S_1,S_2,\ldots,S_{2m}\in C)\nonumber\\
&\geq \alpha\PP(S_2,S_4,\ldots,S_{2m}\in C).
\end{align}
On the other hand, we clearly have
\begin{equation}\label{down}
\PP(T_C>n)\leq\PP(S_2,S_4,\ldots,S_{2m-2}\in C).
\end{equation}
Hence, the problem reduces to that of finding a good estimate of $\PP(S_2,S_4,\ldots, S_{2m}\in C)$; 
the benefit of this reduction is
that the cone $C$ is in general position with respect to the random walk $(S_{2m})_{m\geq 0}$.
One of the key point is that
the increments of $(S_{2m})_{m\geq 0}$ are equal to $0$ with positive probability so that the parity problem does not occur.
The complete proof is not difficult, so we omit it. Now,
the combination of Varopoulos estimate for $\PP(S_2,S_4,\ldots,S_{2m}\in C)$ and inequalities~\eqref{up} and~\eqref{down} gives
the expected result :
$$\gamma^{-1}n^{-\pi/2\beta}\leq\PP(T_C>n)\leq\gamma n^{-\pi/2\beta}$$
for some $\gamma>1$. Therefore, by application of our Theorem~\ref{mainthm} we obtain :

\begin{theo}\label{SR}
Let $(S_n)$ be the simple random walk on $\ZZ^2$ and $C$ be an adapted cone. Then $Q_n\Rightarrow M$.
\end{theo}

\subsection{Degenerated cases}
The question rises whether the sequence of conditional distributions $(Q_n)$ can converge to some limit $Q\not=M$. 
In that case $\PP(T_C>n)$ would not be dominatedly varying.
We do not know any example in dimension 2 (with $C$ adapted to the random walk) but there are some in dimension 3.

\begin{exem} 
Let $(S_n)$ be the simple random walk on $\ZZ^3$, and take $C=\{(x,y,z)\in\RR^3: 0\leq x/2\leq y\leq 2x\}$.
Conditional on $\{T_C>n\}$, the process $\{S_k,k=1\ldots n\}$ is a simple random walk on the axis $\{(0,0,z):z\in\ZZ\}$. Thus $Q_n\Rightarrow Q$, where $Q$ is the law of the process $\{(0,0,\sqrt{3}B_t),t\geq 0\}$, $B_t$ being a stantard one-dimensional Brownian motion. Here $\PP(T_C>n)=(1/3)^n$.
\end{exem}

\begin{exem}
Here again consider the simple random walk on $\ZZ^3$, and take $C=\{(x,y,z)\in\RR^3: 0\leq x/2\leq y\leq 2x\mbox{ and }z\geq 0\}$.
Conditional on $\{T_C>n\}$, the process $\{S_k,k=1\ldots n\}$ is a simple random walk on the axis $\{(0,0,z):z\in\ZZ\}$ conditioned to stay positive. 
Hence $Q_n\Rightarrow Q$, where $Q$ is the law of the process $\{(0,0,\sqrt{3}M_t),t\geq 0\}$, $M_t$ being a one-dimensional Brownian meander. Here $\PP(T_C>n)\sim(1/3)^n(\pi n)^{-1/2}$.
\end{exem}

These examples show that it is possible that $Q_n$ converges to some limit even if $\PP(T_C>n)$ is not dominatedly varying. But then, the limit process is degenerated. The following proposition, which holds in any dimension and for any adapted cone, shows that it is a general fact.

\begin{prop}\label{dgncase}
Suppose $Q_n\Rightarrow Q$ weakly on $\C_1$. If $\PP(T_C>n)$ is not dominatedly varying then 
the limit process lives on the boundary of $C$, i.e.
$$Q(\forall t\in[0,1], w(t)\in \partial C)=1.$$
\end{prop}

The proofs of Theorem~\ref{mainthm} and Proposition~\ref{dgncase} are deferred to Section~\ref{proofs}. Prerequisites are collected in Section~\ref{preparatory}.

\section{Preparatory material}\label{preparatory}

We collect here some of the results  we need to prove Theorem~\ref{mainthm}.

\subsection{More on Brownian meander}

Let $(X_t,0\leq t\leq 1)$ be the canonical process on $\C_1$ for which $X_t(w)=w(t)$.

Given $0<\beta\leq\pi$, let $C\subset\RR^2$ be the cone $\{(r\cos\theta,r\sin\theta):r>0\mbox{ and }0 <\theta<\beta\}$. (One can replace the condition $0 <\theta<\beta$ in the definition of $C$ by $0 \leq\theta<\beta$, $0 <\theta\leq\beta$ or $0 \leq\theta\leq\beta$, for the exit times of these cones are almost surely equal relative to Wiener distribution; hence the distribution of Brownian motion conditional on $\{\tau_C>1\}$ does not depend on that choice.)
Let us denote (as before) $M$ the distribution on $\C_1$ of the Brownian meander of the cone $C$. The Brownian meander is a continuous, non-homogeneous Markov process, with transition density given by
\begin{align}\label{transitions}
M(X_t\in dy)&=e(t,y)dy\\
&= \frac{r^{2\alpha}}{2^{\alpha}\Gamma(\alpha)t^{2\alpha+1}}\exp(-\frac{r^2}{2t})\sin(2\alpha\theta) W^y(\tau_C>1-t)dy\nonumber
\end{align}
for $0< t \leq 1$ and $y=(r\cos\theta,r\sin\theta)$, $r>0$, $0\leq\theta\leq\beta$; $\alpha=\pi/(2\beta)$. And
\begin{align*}
M(X_t\in dy\vert X_s=x)&=p(s,x,t,y)dy\\
&=p^C(t-s,x,y)\frac{W^y(\tau_C>1-t)}{W^x(\tau_C>1-s)}dy\nonumber
\end{align*}
for $0<s<t\leq 1$ and $x,y\in\open{C}$, where $p^C(t,x,y)$ is the heat kernel of $C$.
(See~\cite{Gar09} for a derivation of these formulas; the formula (3.2) given in~\cite{Shi85} for the transition density was misprinted.)
From formula~\eqref{transitions}, one can already obtain a necessary condition for the weak convergence of the conditional distribution $Q_n$ to the Brownian meander.

The conclusion of the following proposition was given by Shimura in~\cite{Shi91} as a consequence of his limit theorem.

\begin{prop}\label{rwexitlaw}
If $Q_n\Rightarrow M$ in $\C_1$, then $\PP(T_C>n)$ is regularly varying with index $\alpha=\pi/(2\beta)$.
\end{prop}
\begin{proof}
Given $t>1$, define $\phi_n(x)=\PP^{x\sqrt{n}}(T_C>[nt]-n)$. By the 
Markov property of the random walk and the definition of $Q_n$, we have
$$\PP(T_C>[nt]\vert T_C>n)=Q_n(\phi_n(X_1)).$$
If $x\in\open{C}$ and $x_n\to x$, then it follows from Donsker's theorem and the Portmanteau theorem (\cite{Bil99}, Theorem 2.1) that
$$\phi_n(x_n)\to \phi(x)=W^x(\tau_C>t-1).$$
Since $Q_n\Rightarrow M$, and $X_1\in\open{C}$ $M$-a.s., the continuous mapping theorem (\cite{Bil99}, Theorem 5.5) shows that
$$Q_n(\phi_n(X_1))\to M(\phi(X_1)).$$
Using~\eqref{transitions}, the limit can be expressed as
$$\int_C e(1,x)W^x(\tau_C>t-1)dx=\int_C te(1,\sqrt{t}y) W^{y}(\tau_C>1-1/t)dy,$$
where we have made the change of variables $x=\sqrt{t}y$ and used the scaling invariance of Brownian motion.
But the last integrand is equal to $t^{-\alpha}e(1/t,y)$, therefore
$$M(\phi(X_1))=t^{-\alpha}.$$
This proves that 
\begin{equation}\label{regalpha}
\PP(T_C>[nt])/\PP(T_C>n)\to t^{-\alpha}
\end{equation}
for all $t>1$. In a similar way, it can be proved that~\eqref{regalpha} also holds for all $t\in(0,1)$, thus $\PP(T_C>n)$ is regularly varying with index $\alpha=\pi/(2\beta)$.
\end{proof}

\subsection{More on conditioned random walk}

Let us finally state here without proof two easy but important facts about the conditioned random walk.
The first one is the Markov property which is inherited from the original unconditioned random walk.
 
Let $\F_t$ be the $\sigma$-algebra generated by the random variables $\{X_s,s\leq t\}$.
The shift operator $\theta_t$ on $\C_1$ is defined by $\theta_t(w)(s)=w(t+s)$.

For all $x\in C$ and $t\in[0,1]$, we denote by $Q_n^{x,t}$ the distribution of $x+\S_n$ conditional on $\{\tau_C(x+\S_n)>t\}$.
Then we have:

\begin{prop}\label{markov}
Let $t=k/n\in[0,1]$ be given. For any $A\in\F_t$ and $B\in \C_1$,
$$Q_n(A;\theta_t^{-1}B)=Q_n(A;Q_n^{X_t,1-t}(B)).$$ 
\end{prop}

The second fact is a limit theorem for the conditioned normalized random walk started inside the cone. 
For $x\in\open{C}$ and $t\in(0,1]$, let $M^{x,t}$ denote the distribution of the standard Brownian motion started at $x$ and conditioned to stay in $C$ until time $t$, that is
$$M^{x,t}(B)=W^x(B\vert \tau_C>t),$$
for any Borel set $B$ of $\C_1$. 
If $t_n=k_n/n\to t$ and $x_n\to x$, then it follows from Donsker's theorem and the Portmanteau theorem that
$$\PP(x_n+\S_n\in B;\tau_C(x_n+\S_n)>t_n)\to W^x(B;\tau_C>t),$$
for any Borel set $B$ such that $\partial B$ is $W^x$-negligible. Since $W^x(\tau_C>t)>0$ (because $x\in\open{C}$), we obtain:

\begin{theo}\label{convbord}
Let $t\in(0,1]$ and $x\in \open{C}$. If $t_n=k_n/n\to t$ and $x_n\to x$, then
$Q_n^{x_n,t_n}\Rightarrow M^{x,t}$.
\end{theo}

Detailed proofs of Proposition~\ref{markov} and Theorem~\ref{convbord} can be found in~\cite{Gar08}.

\section{Proofs of results}\label{proofs}

The proof of Theorem~\ref{mainthm} begins with the analysis of the approach developped by Shimura in~\cite{Shi91}: he proved that the sequence $(Q_n)$ is tight and then observed (without stating it as a general fact) that a sufficient (and necessary) condition for the convergence $Q_n\Rightarrow M$ is that any weak limit point $Q$ of $(Q_n)$ satisfies $Q(w(t)\in\partial C)=0$ for all $t\in(0,1]$. We explain this in Section~\ref{approach}.
We will then show in Section~\ref{unreachability} how this condition relates
to the tails of the exit time. This is the main novelty of this paper.

\subsection{Shimura's approach}\label{approach} 
The paper~\cite{Shi91} of Shimura contains two significant results. 
The first one is the tightness of the sequence $(Q_n)$ under the assumption that the increments of $(S_n)$ are bounded.
Shimura's proof of tightness is quite technical and requires an extension of Theorem~\ref{convbord} to sequences $x_n\to x\in\partial C\setminus\{0\}$.
We do not recall it here and refer the reader to~\cite{Shi91},~Lemma~3.
Let us simply mention that the assumption of bounded increments can be replaced by the following one :
\begin{equation}\label{alternative}
\PP(\max_{i=1\ldots n}\Vert\xi_i\Vert>\sqrt{n}\vert T_C>n)\to 0.
\end{equation}
This is explained in~\cite{Gar08}.
Of course, condition~\eqref{alternative} is satisfied if the steps $\xi_i$ are bounded, but in other cases it is not clear if one should impose an additional moment condition or not.
For example, suppose that we already know that $\PP(T_C>n)\geq \gamma n^{-\alpha}$ for some positive constants $\gamma$ and $\alpha$\footnote{Apart from the random walks with bounded increments that were considered in our Section~\ref{examples}, this lower bound is also known for a certain class of random walks with unbounded increments, see Section~6 of~\cite{Var99}.}. 
A trivial upper bound is given by :
\begin{align*}
\PP(\max_{i=1\ldots n}\Vert\xi_i\Vert>\sqrt{n}\vert T_C>n)&\leq\frac{\PP(\max_{i=1\ldots n}\Vert\xi_i\Vert>\sqrt{n})}{\PP(T_C>n)}\\
&\leq\gamma^{-1}n^{\alpha+1}\PP(\Vert\xi_1\Vert>\sqrt{n}).
\end{align*}
Hence condition~\eqref{alternative} is satisfied whenever $\EE(\Vert\xi_1\Vert^{2\alpha+2})$ is finite. Since $\alpha$ should not be less than $\pi/2\beta$ where $\beta$ is the angle of $C$,
this condition on moments asks for  a finite third moment in the half-plane case (for which we already know that a second moment is sufficient since it is similar to the one-dimensional case), a finite fourth moment in the quarter-plane case, or a finite sixth moment in the octant case. Clearly, we are far from an optimal condition!
In the same spirit, if one knows \textit{a priori} nothing about the asymptotic behavior of $\PP(T_C>n)$, it is still possible to obtain~\eqref{alternative} under a very strong integrability assumption. Namely, the condition is 
\begin{equation}\label{secondalternative}
\EE\left(\Vert\xi_1\Vert^2 a^{\Vert\xi_1\Vert^2}\right)<\infty,
\end{equation}
where $a=\PP(\xi_1\in C)^{-1}$. Indeed, by a classical argument, condition~\eqref{secondalternative} implies that
$na^n\PP(\Vert\xi_1\Vert>\sqrt{n})\to 0$. But, since $C$ is a semi-group, we also have $\PP(T_C>n)\geq a^{-n}$.
Therefore, 
$$\PP(\max_{i=1\ldots n}\Vert\xi_i\Vert>\sqrt{n}\vert T_C>n)
\leq \frac{n\PP(\Vert\xi_1\Vert>\sqrt{n})}{\PP(T_C>n)}\leq na^n\PP(\Vert\xi_1\Vert>\sqrt{n})\to 0,$$
and consequently, the sequence $(Q_n)$ is tight.
We insist on the fact that all these considerations only concern the question of the tightness of $(Q_n)$.

In the rest of this section, we do not assume that the random walk $(S_n)$ has bounded increments, nor any $k$-th moment for $k>2$. 

We now turn to the study of the eventual limit points of $(Q_n)$.
Let $Q$ be a probability measure on $\C_1$. We shall say that $\partial C$ is {\em unreachable} for $Q$ if $Q(w(t)\in\partial C)=0$ for all $t\in(0,1]$.
By virtue of~\eqref{transitions}, the boundary of $C$ is unreachable for $M$. Perhaps, the most striking result in Shimura's paper is that $M$ is the only possible subsequential limit of $(Q_n)$ for which $\partial C$ is unreachable.
The proof is illuminating and we believe it is of interest to reproduce it here.
\begin{prop}[\cite{Shi91}, Proof of Theorem 1]\label{ShiB}
 Let $Q$ be a limit point of the sequence $(Q_n)$.
Then $Q=M$  if and only if $\partial C$ is unreachable for $Q$.  
\end{prop}
\begin{proof}
Let $Q$ be a weak limit point of the sequence $(Q_n)$. 
There exists a subsequence $(Q_{n'})$ which converges weakly to $Q$. To simplify the notation, 
we shall suppose that the whole sequence $(Q_n)$ converges to $Q$.
We will show that $Q$ and $M$ have the same one-dimensional distributions; the generalization to other finite-dimensional distributions is straightforward.

Fix $t\in (0,1]$ and let $f$ be a bounded continuous real function. We have to show that
$$Q(f(X_t))=M(f(X_t)).$$
First choose $0<\lambda< t$ and set $\lambda_n=[n\lambda]/n$ and $t_n=[nt]/n$. Note that $\lambda_n\to\lambda$ and $t_n\to t$ as $n\to \infty$. 

Given a vector $u\in\open{C}$ we set $C_{\epsilon}=\epsilon u +C$ and $\Delta_{\epsilon}=\close{C}\setminus \open{C}_{\epsilon}$.
Note that $\cap_{\epsilon>0}\Delta_{\epsilon}=\partial C$.
For all $\epsilon>0$, define
$$J^n_{\epsilon}=Q_n(X_{\lambda_n}\in C_{\epsilon};f(X_{t_n})).$$
Then,
$$\vert Q_n(f(X_{t_n}))-J^n_{\epsilon}\vert\leq K Q_n(X_{\lambda_n}\in \Delta_{\epsilon}),$$
where $K$ is a bound for $\vert f\vert$. By the continuous mapping theorem $Q_n(X_{\lambda_n}\in dx)\Rightarrow Q(X_{\lambda}\in dx)$.
A standard use of the Portmanteau theorem then shows that 
$$\lim_{\epsilon\to 0}\limsup_{n\to\infty}Q_n(X_{\lambda_n}\in \Delta_{\epsilon})\leq Q(X_{\lambda}\in\partial C)=0\;.$$
In addition, (continuous mapping theorem again) we have
$$\lim_{n\to\infty}Q_n(f(X_{t_n}))=Q(f(X_t)).$$
Therefore,
$$
\lim_{\epsilon\to 0}\limsup_{n\to\infty}\vert Q(f(X_{t}))-J^n_{\epsilon}\vert=0.
$$
Thus, it remains to prove that
$$\lim_{\epsilon\to 0}\lim_{n\to\infty}J^n_{\epsilon}=M(f(X_t)).$$ 
By the Markov property of $Q_n$, we have
\begin{align*}
J^n_{\epsilon}&=Q_n\left(X_{\lambda_n}\in C_{\epsilon};Q_n^{X_{\lambda_n},1-\lambda_n}(f(X_{t_n-\lambda_n}))\right)\\
&=Q_n(X_{\lambda_n}\in C_{\epsilon};\phi_n(X_{\lambda_n})),
\end{align*}
where
$$\phi_n(x)=Q_n^{x,1-\lambda_n}(f(X_{t_n-\lambda_n})),\quad x\in\close{C}.$$
If $x_n\to x\in \open{C}$, then by Theorem~\ref{convbord}, 
$$\lim_{n\to\infty}\phi_n(x_n)=M^{x,1-\lambda}(f(X_{t-\lambda}))=:\phi_{\lambda}(x).$$
Hence, if $w$ is such that $w(\lambda)\notin\partial C_{\epsilon}$, and if $w_n\to w$ uniformly on $[0,1]$, then
\begin{equation}\label{cacestbien}
\lim_{n\to\infty}\indic{C_{\epsilon}}(w_n(\lambda_n))\phi_n(w_n(\lambda_n))=
\indic{C_{\epsilon}}(w(\lambda))\phi(w(\lambda)).
\end{equation}
Let $S$ be the set of all $\epsilon>0$ such that $Q(X_{\lambda}\in\partial C_{\epsilon})=0$. The set $(0,\infty)\setminus S$ is at most countable.
Fix $\epsilon\in S$. By the continuous mapping theorem, it follows from~\eqref{cacestbien} that
$$
\lim_{n\to\infty}J^{n}_{\epsilon}=Q(X_{\lambda}\in{C_{\epsilon}};\phi_{\lambda}(X_{\lambda})).
$$
Now, letting $\epsilon\to 0$ through $S$ gives
$$
\lim_{\epsilon\to 0}\lim_{n\to\infty}J^{n}_{\epsilon}=Q(X_{\lambda}\in \open{C};\phi_{\lambda}(X_{\lambda}))=
Q(\phi_{\lambda}(X_{\lambda})),$$
since $Q(X_{\lambda}\in\partial C)=0$. This last expression does not depend on $\lambda$, so we are going to let $\lambda\to 0$.
To do this, select a sequence $\lambda_n\to 0$. Then, with $Q$-probability one, $w(\lambda_n)\in \open{C}$ for all $n$ (by hypothesis), and $w(\lambda_n)\to w(0)=0$ (by continuity of the paths). Therefore, an easy modification of Theorem~\ref{convmeander} using only the scaling invariance of Brownian motion shows that
$$\phi_{\lambda_n}(w(\lambda_n))=M^{w(\lambda_n),1-\lambda_n}(f(X_{t-\lambda_n}))\to M(f(X_t)).$$
This holds for $Q$-almost all $w$, so by the dominated convergence theorem, we have
$$\lim_{n\to\infty}Q(\phi_{\lambda_n}(X_{\lambda_n}))=M(f(X_t)),$$
which completes the proof.
\end{proof}

\begin{rema}
We point out that Proposition~\ref{ShiB} holds in any dimension provided that the Brownian meander, defined as the weak limit of Brownian motion conditioned to stay in the cone for a unit of time, exists. For example, this holds for cones with a smooth boundary (see~\cite{Gar09} for further details).
\end{rema}

\subsection{Tails of the exit time and unreachability of the boundary}\label{unreachability}
In this section, the dimension is an arbitrary integer $d\geq 1$, and $C$ is an adapted cone of $\RR^d$.
We do not assume that the random walk $(S_n)$ has bounded increments.
It turns out that the unreachability of the boundary of $C$ is closely related to the asymptotic behavior of the tail distribution of the exit time $T_C$.

The first lemma gives a sufficient condition for unreachability.
\begin{lemm}\label{direct} Let $Q$ be a limit point of the sequence $(Q_n)$. 
If $\PP(T_C>n)$ is dominatedly varying, then $\partial C$ is unreachable for $Q$.
\end{lemm}

\begin{proof} Without loss of generality, we shall assume that $\sigma^2=1$.
Choose a vector $u\in\open{C}$ and set $C_{\epsilon}=\epsilon u +C$ and $\Delta_{\epsilon}=\close{C}\setminus \open{C}_{\epsilon}$.
Note that $\cap_{\epsilon>0}\Delta_{\epsilon}=\partial C$. 
Fix $0<s<t\leq 1$ and define 
$$
p(n,\epsilon, R)=\PP(\Vert S_{[ns]}\Vert\leq\sqrt{n}R; S_{[nt]}\in\sqrt{n}\Delta_{\epsilon}\vert T_C>n),\quad n,\epsilon,R>0.
$$
We shall first prove that
\begin{equation}\label{topcontrol}
\forall R>0,\quad\lim_{\epsilon\to 0}\limsup_{n\to\infty}p(n,\epsilon,R)=0.
\end{equation}
Since $\PP(T_C>n)$ is dominatedly varying, there exists a positive constant $\gamma$ such that
$$\forall n\geq 1,\quad \PP(T_C>[ns])\leq \gamma \PP(T_C>n).$$
By the Markov property of the random walk, we have
\begin{eqnarray*}
\lefteqn{\PP(T_C>n;\Vert S_{[ns]}\Vert\leq\sqrt{n}R;S_{[nt]}\in\sqrt{n}\Delta_{\epsilon})}\\
&\leq &\PP(T_C>[ns];\Vert S_{[ns]}\Vert\leq\sqrt{n}R;\PP^{S_{[ns]}}(S_{k_n}\in\sqrt{n}\Delta_{\epsilon}))\\
&\leq &\PP(T_C>[ns])\sup\{\PP^{x}(S_{k_n}\in\sqrt{n}\Delta_{\epsilon})):\Vert x\Vert\leq \sqrt{n}R\}
\end{eqnarray*}
where  $k_n=[nt]-[ns]$. Thus, with our choice of  $\gamma$, we obtain 
\begin{equation}\label{upbound}
p(n,\epsilon,R)\leq \gamma\sup\{\PP(z+S_{k_n}/\sqrt{n}\in\Delta_{\epsilon})):\Vert z\Vert\leq R\}.
\end{equation}
If $z_n\to z$, then the CLT and the Portmanteau theorem imply
$$\limsup_{n\to\infty}\PP(z_n+S_{k_n}/\sqrt{n}\in\Delta_{\epsilon})\leq \N(\Delta_{\epsilon}),$$
where $\N$ is a normal distribution on $\RR^d$. But, as $\epsilon\downarrow 0$, $\Delta_{\epsilon}$ decreases to $\partial C$, a negligible set with respect to Lebesgue measure. Therefore
$$\lim_{\epsilon\to 0}\limsup_{n\to\infty}\PP(z_n+S_{k_n}/\sqrt{n}\in\Delta_{\epsilon})=0.$$
By a compactness argument, the same result holds for the right hand side of~\eqref{upbound}. Thus~\eqref{topcontrol} holds.

Set $t_n=[nt]/n$ and $s_n=[ns]/n$. Then, with regard to $(Q_n)$, relation~\eqref{topcontrol} translates into 
$$\forall R>0,\quad \lim_{\epsilon\to 0}\limsup_{n\to\infty}Q_n(\Vert w(s_n)\Vert\leq R;w(t_n)\in\Delta_{\epsilon})=0.$$
By the Portmanteau theorem, this implies
$$\forall R>0,\quad Q(\Vert w(s)\Vert\leq R;w(t)\in\partial C)=0.$$
Letting $R\to\infty$ completes the proof.
\end{proof}

If the whole sequence $(Q_n)$ converges weakly to some limit $Q$, a converse to Lemma~\ref{direct} holds:

\begin{lemm}\label{converse} Suppose  $Q_n\Rightarrow Q$ in $\C_1$. If $\PP(T_C>n)$ is not dominatedly varying, then $Q(w(t)\in\partial C)=1$ for all $t\in(0,1]$.
\end{lemm}

\begin{proof} Here again, we shall assume without loss of generality that $\sigma^2=1$. Suppose there exists $t\in(0,1]$ such that $Q(w(t)\in\partial C)<1$ and set $t_n=[nt]/n$.
Define $C_{\epsilon}$ and $\Delta_{\epsilon}$ as in the proof of Lemma~\ref{direct}.
Then, by the Portmanteau theorem,
$$\limsup_{n\to\infty}Q_n(w(t_n)\in\Delta_{\epsilon})\leq Q(w(t)\in\Delta_{\epsilon}).$$
Hence, there exists $\epsilon>0$ such that
\begin{equation}\label{petit}
\limsup_{n\to\infty}Q_n(w(t_n)\in\Delta_{\epsilon})<1.
\end{equation}
We shall prove that this implies the dominated variation of $\PP(T_C>n)$.

Let $s\in(0,1)$ be given and set $m=[ns]$. By the Markov property of the random walk, we have
\begin{align}\label{gras}
\PP(T_C>n)&\geq \PP(T_C>n; S_{[mt]}\in\sqrt{m}C_{\epsilon})\nonumber\\
&\geq \PP(T_C>[mt]; S_{[mt]}\in\sqrt{m}C_{\epsilon};\PP^{S_{[mt]}}(T_C>n-[mt]))\nonumber\\
&\geq \PP(T_C>m; S_{[mt]}\in\sqrt{m}C_{\epsilon})\inf_{x\in\sqrt{m}C_{\epsilon}}\PP^{x}(T_C>n-[mt]).
\end{align}
But, since $C$ is a semi-group, 
$$p_m:=\inf_{x\in\sqrt{m}C_{\epsilon}}\PP^{x}(T_C>n-[mt])\geq\PP^{\sqrt{m}\epsilon u}(T_C>n-[mt]),$$
and it follows from Donsker's invariance principle and the Portmanteau theorem  that
\begin{equation}\label{gros}
\liminf_{n\to\infty} p_m\geq p=W^{\epsilon u}(\tau_{\open{C}}>s^{-1}-t)>0.
\end{equation}
Now, dividing both sides of~\eqref{gras} by $\PP(T_C>m)$ gives
$$\frac{\PP(T_C>n)}{\PP(T_C>m)}\geq Q_m(w(t_m)\in C_{\epsilon})\times p_m.$$
Thus, by virtue of~\eqref{petit} and~\eqref{gros}, we obtain
$$\liminf_{n\to\infty}\frac{\PP(T_C>n)}{\PP(T_C>m)}>0.$$
Therefore $\PP(T_C>n)$ is dominatedly varying.
\end{proof}

\subsection{Final steps}
Let us finally give the proofs of our main results, namely Theorem~\ref{mainthm} and Proposition~\ref{dgncase}.

\begin{proof}[Proof of Theorem~\ref{mainthm}] Suppose that the random walk $(S_n)$ has bounded increments. Then, by Lemma~3 in~\cite{Shi91}, we know that the sequence $(Q_n)$ is tight on $\C_1$.
If $\PP(T_C>n)$ is dominatedly varying,  then it follows from Theorem~\ref{ShiB} and Lemma~\ref{direct} that $M$ is the only possible limit point of $(Q_n)$.
Therefore $Q_n\Rightarrow M$.

The converse part is Proposition~\ref{rwexitlaw}, which also gives the index of regular variation $\alpha=\pi/2\beta$.
Note that the dominated variation of $\PP(T_C>n)$ can be derived from Lemma~\ref{converse}, since we know from~\eqref{transitions} that $M(w(t)\in\partial C)=0$ for all $t\in(0,1]$. 
\end{proof}

\begin{proof}[Proof of Proposition~\ref{dgncase}]
Assume the hypotheses of Proposition~\ref{dgncase} are satisfied. Then Lemma~\ref{converse} ensures that
$Q(w(t)\in\partial C)=1$
for all $t\in(0,1]$. Thus
$$Q(\forall t\in\QQ\cap [0,1], w(t)\in\partial C)=1.$$
The result follows since every path $w\in\C_1$ is continuous, $\QQ$ is dense in $[0,1]$ and $\partial C$ is a closed set.
\end{proof}

\section*{Acknowledgments}
I would like to thank E.~Lesigne, M.~Peign\'e and an anonymous referee for helping in many ways to improve the original manuscript.

\end{document}